\documentclass[a4paper,11pt]{amsart}

\usepackage{graphicx}
\usepackage{mathptmx}
\usepackage{amsmath}
\usepackage{amssymb}
\usepackage{enumitem}
\usepackage{xcolor}

\newmuskip\pFqmuskip

\newcommand*\pFq[6][8]{%
  \begingroup 
  \pFqmuskip=#1mu\relax
  \mathcode`=\string"8000
  \begingroup\lccode`\~=`\,
  \lowercase{\endgroup\let~}\pFqcomma
  F^{#2}_{#3}{\left(\genfrac..{0pt}{}{#4}{#5}\bigg|#6\right)}%
  \endgroup
}
\newcommand{\pFqcomma}{\mskip\pFqmuskip}

\newtheorem{theorem}{Theorem}[section]

\begin{document}

\title[Generalization of Spivey's recurrence relation]{Generalization of Spivey's recurrence relation }

\author{Taekyun  Kim}
\address{Department of Mathematics, Kwangwoon University, Seoul 139-701, Republic of Korea}
\email{tkkim@kw.ac.kr}
\author{Dae San  Kim}
\address{Department of Mathematics, Sogang University, Seoul 121-742, Republic of Korea}
\email{dskim@sogang.ac.kr}


\maketitle

\begin{abstract}
In 2008, Spivey found a recurrence relation for the Bell numbers $\phi_{n}$. We consider the probabilistic $r$-Bell polynomials associated with $Y$, $\phi_{n,r}^{Y}(x)$, which are a probabilistic extension of the $r$-Bell polynomials. Here $Y$ is a random variable whose moment generating function exists in some neighborhood of the origin and $\phi_{n}=\phi_{n,0}^{1}(1)$. The aim of this paper is to generalize the relation for the Bell numbers to that for the probabilistic $r$-Bell polynomials associated with $Y$.
\end{abstract}

\section{Introduction}
The Stirling number ${n \brace k}$ counts the number of ways to partition a set of $n$ objects into $k$ nonempty subsets. The Bell polynomials are defined by $\phi_{n}(x)=\sum_{k=0}^{n}{n \brace k}x^{k}$, with the Bell numbers given by $\phi_{n}=\phi_{n}(1)$. For $r\in\mathbb{N}$, the $r$-Stirling number ${n+r \brace k+r}_{r}$ enumerates the number of ways to partition  a set of $n+r$ objects into $k+r$ nonempty subsets such that the first $r$ elements are in distinct subsets.
The $r$-Bell polynomials are defined by $\phi_{n,r}(x)=\sum_{k=0}^{n}{n+r \brace k+r}_{r}x^{k}$, with the $r$-Bell numbers given by $\phi_{n,r}=\phi_{n,r}(1)$. Assume that $Y$ is a random variable whose moment generating function exists in some neighbourhood of the origin. Let $(Y_{j})_{j\ge 1}$ be a sequence of mutually independent copies of the random variable $Y$, and let $S_{k}=Y_{1}+Y_{2}+\cdots+Y_{k},\ (k\ge 1)$, with $S_{0}=0$. Then we define the probabilistic $r$-Stirling numbers of the second kind associated with $Y$, ${n+r \brace k+r}_{r,Y}$, as a probabilistic extension of the $r$-Stirling numbers of the second kind ${n+r \brace k+r}_{r}$. The probabilistic $r$-Bell polynomials associated with $Y$ are defined by $\phi_{n,r}^{Y}(x)=\sum_{k=0}^{n}{n+r \brace k+r}_{r,Y}x^{k}$, with
the probabilistic $r$-Bell numbers associated with $Y$ given by $\phi_{n,r}^{Y}=\phi_{n,r}^{Y}(1)$. \par
Spivey found the following recurrence relation on the Bell numbers:
\begin{equation*}
\phi_{n+k}=\sum_{l=0}^{k}\sum_{j=0}^{n}j^{k-l}\binom{k}{l}{n \brace j}\phi_{l}.
\end{equation*}
The aim of this paper is to generalize the relation for the Bell numbers to that for the probabilistic $r$-Bell polynomials associated with $Y$. Indeed, we obtain
\begin{equation}
\phi_{j+n,r}^{Y}(y)=\sum_{l=0}^{n}\binom{n}{l}r^{n-l}\sum_{k=0}^{l}\sum_{m=0}^{j}\frac{y^{k}}{k!}\binom{j}{m}\sum_{l_{1}+\cdots+l_{k}=l}\binom{l}{l_{1},\dots,l_{k}}E\Big[S_{k}^{j-m}\prod_{i=1}^{k}Y_{i}^{l_{i}}\Big]\phi_{m,r}^{Y}(y), \label{0}
\end{equation}
where $l_{1},l_{2},\dots l_{k}$ are positive integers. \\
By using the simple fact \eqref{15}, we show \eqref{0}, first for $\phi_{n}^{Y}=\phi_{n,0}^{Y}(1)$, second for $\phi_{n}^{Y}(y)=\phi_{n,0}^{Y}(y)$, and finally for $\phi_{n,r}^{Y}(y)$. We note here that $\phi_{n}=\phi_{n,0}^{1}(1)$. For the rest of this section, we recall the facts that are needed throughout this paper. \par

\vspace{0.1in}

For $n\ge 0$, the Stirling numbers of the second kind are defined by
\begin{equation}
x^{n}=\sum_{k=0}^{n}{n \brace k}(x)_{k},\quad (\mathrm{see}\ [1-26]), \label{1}	
\end{equation}
where $(x)_{0}=1,\ (x)_{k}=x(x-1)\cdots (x-k+1),\ (k\ge 1)$. \par
The Bell numbers are given by
\begin{equation}
\phi_{n}=\sum_{k=0}^{n}{n\brace k},\quad (n\ge 0),\quad (\mathrm{see}\ [15,21]).\label{2}
\end{equation}
The Bell polynomials are defined by
\begin{equation}
e^{x(e^{t}-1)}=\sum_{n=0}^{\infty}\phi_{n}(x)\frac{t^{n}}{n!}. \label{3}	
\end{equation}
From \eqref{3}, we note that
\begin{equation}
\phi_{n}(x)=\sum_{k=0}^{n}{n\brace k}x^{k},\quad (n\ge 0),\quad (\mathrm{see}\ [7,8,10,15,21]).\label{4}
\end{equation}
Note that, for $x=1$, $\phi_{n}=\phi_{n}(1)$ are the Bell numbers. \par
For $r\in\mathbb{N}$ and $n\ge 0$, the $r$-Stirling numbers of the second kind are given by
\begin{equation}
(x+r)^{n}=\sum_{k=0}^{n}{n+r \brace k+r}_{r}(x)_{k},\quad (\mathrm{see}\ [20]).\label{5}	
\end{equation}
From \eqref{5}, we note that
\begin{align}
\sum_{l=0}^{n}{n+r \brace l+r}_{r}(x)_{l}&=(x+r)^{n}=\sum_{k=0}^{n}\binom{n}{k}r^{n-k}x^{k}\label{6}\\
&=\sum_{k=0}^{n}\binom{n}{k}r^{n-k}\sum_{l=0}^{k}{k \brace l}(x)_{l}\nonumber \\
&=\sum_{l=0}^{n}\sum_{k=l}^{n}\binom{n}{k}{k\brace l}r^{n-k}(x)_{l}.\nonumber
\end{align}
Thus, by comparing the coefficients on both sides of \eqref{6}, we get
\begin{equation}
{n+r \brace l+r}_{r}=\sum_{k=l}^{n}\binom{n}{k}{k \brace l}r^{n-k},\quad (n\ge l\ge 0).\label{7}
\end{equation}
It is known that the Bell numbers satisfy the recurrence relation
\begin{equation}
\phi_{n+1}=\sum_{k=0}^{n}\binom{n}{k}\phi_{k},\quad (n\ge 0),\quad (\mathrm{see}\ [14,18]).\label{8}
\end{equation}
In 2008, Spivey discovered the following remarkable relation
\begin{equation}
\phi_{n+k}=\sum_{l=0}^{k}\sum_{j=0}^{n}j^{k-l}\binom{k}{l}{n \brace j}\phi_{l},\quad (\mathrm{see}\ [24]). \label{9}
\end{equation}
The $r$-Bell numbers are defined by
\begin{equation}
\phi_{n,r}=\sum_{k=0}^{n}{n+r \brace k+r}_{r},\quad (n\ge 0),\quad (\mathrm{see}\ [13-20]).\label{10}
\end{equation}
The $r$-Bell polynomials are given by
\begin{equation}
e^{rt}e^{x(e^{t}-1)}=\sum_{n=0}^{\infty}\phi_{n,r}(x)\frac{t^{n}}{n!},\quad (\mathrm{see}\ [13-20]).\label{11}
\end{equation}
Thus, by \eqref{11}, we get
\begin{equation}
\phi_{n,r}(x)=\sum_{k=0}^{n}{n+r \brace k+r}_{r}x^{k},\quad (n\ge 0),\quad (\mathrm{see}\ [6,8,16,17]).\label{12}
\end{equation}
Note that, for $x=1$, $\phi_{n,r}=\phi_{n,r}(1)$ are the $r$-Bell numbers. \par
For any integer $k\ge 0$, the partial Bell polynomials are given by
\begin{equation}
\frac{1}{k!}\bigg(\sum_{m=1}^{\infty}x_{m}\frac{t^{n}}{n!}\bigg)^{k}=\sum_{n=k}^{\infty}B_{n,k}\big(x_{1},x_{2},\dots,x_{n-k+1}\big)\frac{t^{n}}{n!},\quad (\mathrm{see}\ [10,11,14]),\label{13}	
\end{equation}
where
\begin{equation}
\begin{aligned}
&B_{n,k}(x_{1},x_{2},\dots,x_{n-k+1})\\
&=\sum_{\substack{l_{1}+\cdots+l_{n-k+1}\\ l_{1}+2l_{2}+\cdots+(n-k+1)l_{n-k+1}=n}}\frac{n!}{l_{1}!l_{2}!\cdots l_{n-k+1}!}\bigg(\frac{x_{1}}{1!}\bigg)^{l_{1}} \bigg(\frac{x_{2}}{2!}\bigg)^{l_{2}}\cdots \bigg(\frac{x_{n-k+1}}{(n-k+1)!}\bigg)^{l_{n-k+1}}.
\end{aligned}	\label{14}
\end{equation}
We will use the following simple fact
\begin{align}
f(x+t)&=\sum_{n=0}^{\infty}\frac{f^{(n)}(x)}{n!}(x+t-x)^{n}=\sum_{n=0}^{\infty}\frac{f^{(n)}(x)}{n!}t^{n}\label{15}\\
&=\sum_{n=0}^{\infty}\frac{t^{n}D_{x}^{n}}{n!}f(x)=e^{tD_{x}}f(x), \nonumber	
\end{align}
where $D_{x}=\frac{d}{dx}$. \par
Assume that $Y$ is a random variable such that the moment generating function of $Y$,
\begin{equation*}
E\big[e^{tY}\big]=\sum_{n=0}^{\infty}E[Y^{n}]\frac{t^{n}}{n!},\quad (|t|<r),
\end{equation*}
exists for some $r>0$ (see [3,14,22]). Here $E$ denotes the mathematical expectation. Let $(Y_{j})_{j\ge 1}$ be a sequence of mutually independent copies of the random variable $Y$, and let $S_{k}=Y_{1}+Y_{2}+\cdots+Y_{k},\ (k\ge 1)$, with $S_{0}=0$. \par
Then we recall the probabilistic Stirling numbers of the second kind associated with $Y$, ${n\brace k}_{Y}$, which are defined in terms of the $n$-th moments of $S_{l},\ (l=0,1,2,\dots,k)$:
\begin{equation}
{n \brace k}_{Y}=\frac{1}{k!}\sum_{l=0}^{k}\binom{k}{l}(-1)^{k-l}E\big[S_{l}^{n}\big],\quad (0\le k\le n),\quad (\mathrm{see}\ [3,15]). \label{16}
\end{equation}
The probabilistic Bell polynomials associated with $Y$ are defined by
\begin{equation}
\phi_{n}^{Y}(x)=\sum_{k=0}^{n}{n\brace k}_{Y}x^{k},\quad (n\ge 0),\quad (\mathrm{see}\ [14,23]). \label{17}
\end{equation}
From [14], we note that the probabilistic Bell polynomials associated with $Y$ satisfy the recurrence relation
\begin{equation}
\phi_{n+1}^{Y}(x)=x\sum_{k=0}^{n}\binom{n}{k}E\big[Y^{k+1}\big]\phi_{n-k}^{Y}(x).\label{18}
\end{equation}

\section{Generalization of Spivey's recurrence relation}
Let $(Y_{j})_{j\ge 1}$ be a sequence of mutually independent copies of the random variable $Y$, and let
\begin{equation*}
S_{0}=0,\quad S_{k}=Y_{1}+Y_{2}+\cdots+Y_{k},\quad (k\in\mathbb{N}). 	
\end{equation*}
For $r\in\mathbb{N}$, we define the {\it{probabilistic $r$-Stirling numbers of the second kind associated with $Y$}} by
\begin{equation}
\frac{1}{k!}\Big(E[e^{tY}]-1\Big)^{k}e^{rt}=\sum_{n=k}^{\infty}{n+r\brace k+r}_{r,Y}\frac{t^{n}}{n!},\quad (k\ge 0). \label{19}	
\end{equation}
When $Y=1$, ${n+r \brace k+r}_{r,Y}={n+r \brace k+r}_{r},\ (n\ge k\ge 0)$. \par
From \eqref{19}, we note that
\begin{align}
\sum_{n=k}^{\infty}{n+r\brace k+r}_{r,Y}\frac{t^{n}}{n!}&=\frac{1}{k!}\Big(E\big[e^{Yt}\big]-1\Big)^{k}e^{rt}\label{20}\\
&=\frac{1}{k!}\sum_{j=0}^{k}\binom{k}{j}(-1)^{k-j}E\big[e^{(Y_{1}+\cdots+Y_{j})t}\big]e^{rt}\nonumber \\
&=\frac{1}{k!}\sum_{j=0}^{k}\binom{k}{j}(-1)^{k-j}E\big[e^{(S_{j}+r)t}\big] \nonumber \\
&=\sum_{n=0}^{\infty}\bigg(\frac{1}{k!}\sum_{j=0}^{k}\binom{k}{j}(-1)^{k-j}E\big[(S_{j}+r)^{n}\big]\bigg)\frac{t^{n}}{n!}. \nonumber	
\end{align}
Therefore, by comparing the coefficients on both sides of \eqref{20}, we obtain the following theorem.
\begin{theorem}
For $n\ge k\ge 0$, we have
\begin{displaymath}
{n+r\brace k+r}_{r,Y}= \frac{1}{k!}\sum_{j=0}^{k}\binom{k}{j}(-1)^{k-j}E\big[(S_{j}+r)^{n}\big].
\end{displaymath}
\end{theorem}
In view of \eqref{12}, we define the {\it{probabilistic $r$-Bell polynomials associated with $Y$}} by
\begin{equation}
\phi_{n,r}^{Y}(x)=\sum_{k=0}^{n}{n+r \brace k+r}_{r,Y}x^{k},\quad (n\ge 0). \label{21}
\end{equation}
When $Y=1$, $\phi_{n,r}^{Y}(x)=\phi_{n,r}(x)$. For $x=1$, $\phi_{n,r}^{Y}=\phi_{n,r}^{Y}(1)$ are called the probabilistic $r$-Bell numbers associated with $Y$. \par
From \eqref{21}, we have
\begin{align}
e^{x(E[e^{Yt}]-1)} e^{rt}&=\sum_{k=0}^{\infty}\frac{x^{k}}{k!}\Big(E[e^{Yt}]-1\Big)^{k}e^{rt}\label{22} \\
&=\sum_{k=0}^{\infty}x^{k}\sum_{n=k}^{\infty}{n+r\brace k+r}_{r,Y}\frac{t^{n}}{n!}=\sum_{n=0}^{\infty}\sum_{k=0}^{n}x^{k}{n+r \brace k+r}_{r,Y}\frac{t^{n}}{n!}\nonumber \\
&=\sum_{n=0}^{\infty}\phi_{n,r}^{Y}(x)\frac{t^{n}}{n!}.\nonumber
\end{align}
Therefore, by \eqref{22}, we obtain the following theorem.
\begin{theorem}
The generating function of the probabilistic $r$-Bell polynomials associated with $Y$ is given by
\begin{equation}
e^{x(E[e^{Yt}]-1)} e^{rt}= \sum_{n=0}^{\infty}\phi_{n,r}^{Y}(x)\frac{t^{n}}{n!}.\label{23}
\end{equation}
\end{theorem}
From \eqref{23}, we have
\begin{align}
&\sum_{n=0}^{\infty}\phi_{n,r}^{Y}(x)\frac{t^{n}}{n!}=e^{x(E[e^{Yt}]-1)} e^{rt}=\sum_{k=0}^{\infty}\frac{1}{k!}\bigg(\sum_{j=1}^{\infty}xE[Y^{j}]\frac{t^{j}}{j!}\bigg)^{k}e^{rt}\label{24}\\
&=\sum_{l=0}^{\infty}\sum_{k=0}^{l}B_{l,k}\Big(xE[Y],xE[Y^{2}],\dots,xE[Y^{l-k+1}]\Big)\frac{t^{l}}{l!}\sum_{m=0}^{\infty}r^{m}\frac{t^{m}}{m!}\nonumber \\
&=\sum_{n=0}^{\infty}\sum_{l=0}^{n}\binom{n}{l} \sum_{k=0}^{l}B_{l,k}\Big(xE[Y],xE[Y^{2}],\dots,xE[Y^{l-k+1}]\Big)r^{n-l}\frac{t^{n}}{n!}.\nonumber
\end{align}
Therefore, by comparing the coefficients on both sides of \eqref{24}, we obtain the following theorem.
\begin{theorem}
For $n\ge 0$, we have
\begin{displaymath}
\phi_{n,r}^{Y}(x)=\sum_{l=0}^{n}\binom{n}{l} \sum_{k=0}^{l}B_{l,k}\Big(xE[Y],xE[Y^{2}],\dots,xE[Y^{l-k+1}]\Big)r^{n-l}.
\end{displaymath}
\end{theorem}
By \eqref{23}, we get
\begin{align}
\sum_{n=0}^{\infty}\phi_{n+1,r}^{Y}(x)\frac{t^{n}}{n!}&=\frac{d}{dt}\sum_{n=0}^{\infty}\phi_{n,r}^{Y}(x)\frac{t^{n}}{n!}=\frac{d}{dt}\bigg(e^{x(E[e^{tY}]-1)} e^{rt}\bigg) \label{25} \\
&=xE[Ye^{Yt}]e^{x(E[e^{Yt}]-1)}e^{rt}+r e^{x(E[e^{tY}]-1)}e^{rt} \nonumber \\
&=x\sum_{k=0}^{\infty}E[Y^{k+1}]\frac{t^{k}}{k!}\sum_{m=0}^{\infty}\phi_{m,r}^{Y}(x)\frac{t^{m}}{m!}+r\sum_{n=0}^{\infty}\phi_{n,r}^{Y}(x)\frac{t^{n}}{n!} \nonumber \\
&=\sum_{n=0}^{\infty}\bigg(x\sum_{k=0}^{n}\binom{n}{k}E[Y^{k+1}]\phi_{n-k,r}^{Y}(x)+r\phi_{n,r}^{Y}(x)\bigg)\frac{t^{n}}{n!}.\nonumber
\end{align}
Therefore, by comparing the coefficients on both sides of \eqref{25}, we obtain the following theorem.
\begin{theorem}
For $n\ge 0$, we have
\begin{displaymath}
\phi_{n+1,r}^{Y}(x)= x\sum_{k=0}^{n}\binom{n}{k}E[Y^{k+1}]\phi_{n-k,r}^{Y}(x)+r\phi_{n,r}^{Y}(x).
\end{displaymath}
\end{theorem}
From \eqref{17}, we note that
\begin{equation}
e^{E[e^{Yx}]-1}=\sum_{n=0}^{\infty}\phi_{n}^{Y}\frac{x^{n}}{n!}. \label{26}	
\end{equation}
By using \eqref{15}, we have
\begin{align}
e^{tD_{x}}e^{E[e^{Yx}]-1}&=e^{E[e^{Y(x+t)}]-1}\label{27} \\
&=e^{E[e^{Yx}(e^{Yt}-1)+e^{Yx}] -1}\nonumber \\
 &=e^{E[e^{Yx}(e^{Yt}-1)]}e^{E[e^{Yx}]-1}.\nonumber
\end{align}
Now, we observe that
\begin{align}
&e^{E[e^{Yx}(e^{Yt}-1)]}=\sum_{k=0}^{\infty}\frac{1}{k!}\Big(E\big[e^{Yx} (e^{Yt}-1)\big]\Big)^{k} \label{27-1} \\
&=\sum_{k=0}^{\infty}\frac{1}{k!}E\Big[e^{(Y_{1}+\cdots+Y_{k})x}(e^{Y_{1}t}-1)(e^{Y_{2}t}-1)\cdots (e^{Y_{k}t}-1)\Big]\nonumber \\
&=\sum_{k=0}^{\infty}\frac{1}{k!}\sum_{n=k}^{\infty}\sum_{l_{1}+\cdots+l_{k}=n}\binom{n}{l_{1},\dots,l_{k}}E\big[e^{S_{k}x}Y_{1}^{l_{1}}Y_{2}^{l_{2}}\cdots Y_{k}^{l_{k}}\big]\frac{t^{n}}{n!}\nonumber \\
&=\sum_{n=0}^{\infty}\sum_{k=0}^{n}\frac{1}{k!}\sum_{l_{1}+\cdots+l_{k}=n}\binom{n}{l_{1},\dots,l_{k}}E\bigg[e^{S_{k}x}\prod_{i=1}^{k}Y_{i}^{l_{i}}\bigg]\frac{t^{n}}{n!},\nonumber
\end{align}
where $l_{1},l_{2},\dots l_{k}$ are positive integers.
From \eqref{26}, \eqref{27} and \eqref{27-1}, we note that
\begin{align}
&e^{tD_{x}}e^{E[e^{Yx}]-1}= e^{E[e^{Yx}(e^{Yt}-1)]}e^{E[e^{Yx}]-1} \label{28} \\
&=\sum_{n=0}^{\infty}\sum_{k=0}^{n}\frac{1}{k!}\sum_{l_{1}+\cdots+l_{k}=n}\binom{n}{l_{1},\dots,l_{k}}\frac{t^{n}}{n!}\sum_{j=0}^{\infty}E\bigg[S_{k}^{j}\prod_{i=1}^{k}Y_{i}^{l_{i}}\bigg]\frac{x^{j}}{j!}\sum_{m=0}^{\infty}\phi_{m}^{Y}\frac{x^{m}}{m!}\nonumber \\
&=\sum_{n=0}^{\infty}\sum_{l=0}^{\infty}\sum_{k=0}^{n}\sum_{m=0}^{l}\binom{l}{m}\frac{1}{k!}\sum_{l_{1}+\cdots+l_{k}=n}\binom{n}{l_{1},\dots,l_{k}}E\bigg[S_{k}^{l-m}\prod_{i=1}^{k}Y_{i}^{l_{i}}\bigg]\phi_{m}^{Y}\frac{t^{n}}{n!}\frac{x^{l}}{l!},\nonumber
\end{align}
where $l_{1},l_{2},\dots l_{k}$ are positive integers.
On the other hand, by \eqref{26}, we get
\begin{align}
&e^{tD_{x}}e^{E[e^{Yx}]-1}=e^{tD_{x}}\sum_{l=0}^{\infty}\phi_{l}^{Y}\frac{x^{l}}{l!}=\sum_{n=0}^{\infty}\frac{t^{n}}{n!}D_{x}^{n}\sum_{l=0}^{\infty}\phi_{l}^{Y}\frac{x^{l}}{l!}\label{29} \\
&=\sum_{n=0}^{\infty}\sum_{l=0}^{\infty}\phi_{l+n}^{Y}\frac{t^{n}}{n!}\frac{x^{l}}{l!}.\nonumber
\end{align}
Therefore, by \eqref{28} and \eqref{29}, we obtain the following theorem.
\begin{theorem}
For $n,l\ge 0$, we have
\begin{displaymath}
\phi_{l+n}^{Y}=\sum_{k=0}^{n}\sum_{m=0}^{l}\binom{l}{m}\frac{1}{k!}\sum_{l_{1}+\cdots+l_{k}=n}\binom{n}{l_{1},\dots,l_{k}} E\bigg[S_{k}^{l-m}\prod_{i=1}^{k}Y_{i}^{l_{i}}\bigg]\phi_{m}^{Y},
\end{displaymath}
where $l_{1},l_{2},\dots l_{k}$ are positive integers.
\end{theorem}
We note here that this reduces to the  Spivey's relation \eqref{9} if $Y=1$.
\begin{align*}
\phi_{l+n}&=\sum_{k=0}^{n}\sum_{m=0}^{l}\binom{l}{m}\frac{1}{k!}\sum_{l_{1}+\cdots+l_{k}=n}\binom{n}{l_{1},\dots,l_{k}}k^{l-m}\phi_{m} \\
&=\sum_{k=0}^{n}\sum_{m=0}^{l}\binom{l}{m}{n \brace k}k^{l-m}\phi_{m}.
\end{align*}
We extend Theorem 2.5 to the probabilistic Bell polynomials associated with $Y$, which are given by (see \eqref{17}, Theorem 2.2)
\begin{equation}
e^{x(E[e^{Yt}]-1)}=\sum_{n=0}^{\infty}\phi_{n}^{Y}(x)\frac{t^{n}}{n!}.\label{30}	
\end{equation}
Observe that
\begin{align}
e^{tD_{x}} e^{y(E[e^{Yx}]-1)}&=e^{y(E[e^{Y(x+t)}]-1)}\label{31}	\\
&=e^{y(E[e^{Yx}(e^{Yt}-1)+E[e^{Yx}] -1)}\nonumber\\
&=e^{yE[e^{Yx}(e^{Yt}-1)]}e^{y(E[e^{Yx}]-1)}.\nonumber
\end{align}
Now, we note that
\begin{align}
&e^{yE[e^{Yx}(e^{Yt}-1)]}=\sum_{k=0}^{\infty}\frac{y^{k}}{k!}\Big(E\big[e^{Yx} (e^{Yt}-1)\big]\Big)^{k}\label{32}	\\
&=\sum_{k=0}^{\infty}\frac{y^{k}}{k!}E\Big[e^{(Y_{1}+Y_{2}+\cdots+Y_{k})x} (e^{Y_{1}t}-1) (e^{Y_{2}t}-1)\cdots (e^{Y_{k}t}-1)\Big] \nonumber \\
&= \sum_{k=0}^{\infty}\frac{y^{k}}{k!}E\Big[e^{S_{k}x}(e^{Y_{1}t}-1) (e^{Y_{2}t}-1)\cdots (e^{Y_{k}t}-1)\Big]\nonumber \\
&= \sum_{k=0}^{\infty}\frac{y^{k}}{k!}\sum_{n=k}^{\infty}\sum_{l_{1}+\cdots+l_{k}=n}\binom{n}{l_{1},\dots,l_{k}}E\Big[e^{S_{k}x}Y_{1}^{l_{1}}Y_{2}^{l_{2}}\cdots Y_{k}^{l_{k}}\Big]\frac{t^{n}}{n!}\nonumber \\
&=\sum_{n=0}^{\infty}\bigg(\sum_{k=0}^{n}\frac{y^{k}}{k!} \sum_{l_{1}+\cdots+l_{k}=n}\binom{n}{l_{1},\dots,l_{k}}E\bigg[e^{S_{k}x}\prod_{i=1}^{k}Y_{i}^{l_{i}}\bigg]\frac{t^{n}}{n!},\nonumber
\end{align}
where $l_{1},l_{2},\dots l_{k}$ are positive integers.
From \eqref{31} and \eqref{32}, we have
\begin{align}
&e^{tD_{x}}e^{y(E[e^{Yx}]-1)}= e^{yE[e^{Yx}(e^{Yt}-1)]} e^{y(E[e^{Yx}]-1)} \label{33}\\
&=\sum_{n=0}^{\infty}\sum_{k=0}^{n}\frac{y^{k}}{k!}\sum_{l_{1}+\cdots+l_{k}=n}\binom{n}{l_{1},\dots,l_{k}}\frac{t^{n}}{n!}\sum_{j=0}^{\infty}E\bigg[S_{k}^{j}\prod_{i=1}^{k}Y_{i}^{l_{i}}\bigg]\frac{x^{j}}{j!}\sum_{m=0}^{\infty}\phi_{m}^{Y}(y)\frac{x^{m}}{m!}\nonumber \\
&=\sum_{n=0}^{\infty}\sum_{l=0}^{\infty}\sum_{k=0}^{n}\sum_{m=0}^{l} \frac{y^{k}}{k!}\sum_{l_{1}+\cdots+l_{k}=n}\binom{n}{l_{1},\dots,l_{k}}\binom{l}{m} E\bigg[S_{k}^{l-m}\prod_{i=1}^{k}Y_{i}^{l_{i}}\bigg]\phi_{m}^{Y}(y)\frac{t^{n}}{n!}\frac{x^{l}}{l!},\nonumber
\end{align}
where $l_{1},l_{2},\dots l_{k}$ are positive integers.
On the other hand, by \eqref{30}, we get
\begin{align}
e^{tD_{x}}e^{y(E[e^{Yx}]-1)}&=\sum_{n=0}^{\infty}\frac{D_{x}^{n}}{n!}t^{n}\sum_{l=0}^{\infty}\phi_{l}^{Y}(y)\frac{x^{l}}{l!}\label{34} \\
&=\sum_{n=0}^{\infty}\sum_{l=0}^{\infty}\phi_{l+n}^{Y}(y)\frac{t^{n}}{n!}\frac{x^{l}}{l!}.\nonumber
\end{align}
Therefore, by \eqref{33} and \eqref{34}, we obtain the following theorem.
\begin{theorem}
For $n,l\ge 0$, we have
\begin{displaymath}
\phi_{l+y}^{Y}(y)=\sum_{k=0}^{n}\sum_{m=0}^{l}\binom{l}{m}\frac{y^{k}}{k!} \sum_{l_{1}+\cdots+l_{k}=n}\binom{n}{l_{1},\dots,l_{k}}E\bigg[S_{k}^{l-m}\prod_{i=1}^{k}Y_{i}^{l_{i}}\bigg]\phi_{m}^{Y}(y),
\end{displaymath}
where $l_{1},l_{2},\dots l_{k}$ are positive integers.
\end{theorem}
When $Y=1$, we note that
\begin{align}
\phi_{l+n}(y)&=\sum_{k=0}^{n}\sum_{m=0}^{l}\binom{l}{m}y^{k}\frac{1}{k!}\sum_{l_{1}+\cdots+l_{k}=n}\binom{n}{l_{1},\dots,l_{k}}k^{l-m}\phi_{m}(y)\label{34-1}\\
&=\sum_{k=0}^{n}\sum_{m=0}^{l}\binom{l}{m}y^{k}{n \brace k}k^{l-m}\phi_{m}(y), \nonumber
\end{align}
where $l,n$ are nonnegative integers. Note that \eqref{34-1} reduces to the Spivey's relation \eqref{9} by letting $y=1$. \par
From \eqref{15} and \eqref{23}, we note that
\begin{align}
&e^{tD_{x}}\Big(e^{y(E[e^{Yx}]-1)}e^{rx}\Big)= e^{y(E[e^{Y(x+t)}]-1)}e^{r(x+t)} \label{35} \\
&= e^{y(E[e^{Yx}(e^{Yt}-1)])} e^{rt}e^{y(E[e^{Yx}]-1)} e^{rx} \nonumber \\
&=\sum_{k=0}^{\infty}\frac{y^{k}}{k!} E\Big[e^{(Y_{1}+Y_{2}+\cdots+Y_{k})x} (e^{Y_{1}t}-1) \cdots (e^{Y_{k}t}-1)\Big]e^{rt}  e^{y(E[e^{Yx}]-1)}e^{rx} \nonumber \\
&=\sum_{k=0}^{\infty}\frac{y^{k}}{k!}\sum_{l=k}^{\infty}\sum_{l_{1}+\cdots+l_{k}=l}\binom{l}{l_{1},\dots,l_{k}}E\Big[e^{S_{k}x}\prod_{i=1}^{k}Y_{i}^{l_{i}}\Big]\frac{t^{l}}{l!}e^{rt}e^{y(E[e^{Yx}]-1)} e^{rx} \nonumber\\
&=\sum_{l=0}^{\infty}\sum_{k=0}^{l}\frac{y^{k}}{k!}\sum_{l_{1}+\cdots+l_{k}=l}\binom{l}{l_{1},\dots,l_{k}}\frac{t^{l}}{l!}e^{rt}E\Big[e^{S_{k}x}\prod_{i=1}^{k}Y_{i}^{l_{i}}\Big]e^{y(E[e^{Yx}]-1)} e^{rx} \nonumber\\
&=\sum_{n=0}^{\infty}\sum_{l=0}^{n}\binom{n}{l}r^{n-l}\sum_{k=0}^{l}\frac{y^{k}}{k!}\sum_{l_{1}+\cdots+l_{k}=l}\binom{l}{l_{1},\dots,l_{k}}\frac{t^{n}}{n!}\sum_{j=0}^{\infty}E\Big[S_{k}^{j}\prod_{i=1}^{k}Y_{i}^{l_{i}}\Big]\frac{x^{j}}{j!}\sum_{m=0}^{\infty}\phi_{m,r}^{Y}(y)\frac{x^{m}}{m!} \nonumber \\
&=\sum_{n=0}^{\infty}\sum_{j=0}^{\infty}\sum_{m=0}^{j}\sum_{l=0}^{n}\binom{n}{l}r^{n-l}\binom{j}{m}\sum_{k=0}^{l}\frac{y^{k}}{k!}\sum_{l_{1}+\cdots+l_{k}=l}\binom{l}{l_{1},\dots,l_{k}} E\Big[S_{k}^{j-m}\prod_{i=1}^{k}Y_{i}^{l_{i}}\Big]\phi_{m,r}^{Y}(y)\frac{t^{n}}{n!}\frac{x^{j}}{j!},\nonumber
\end{align}
where $l_{1},l_{2},\dots l_{k}$ are positive integers.
On the other hand, by \eqref{23}, we get
\begin{align}
	e^{tD_{x}}\Big(e^{y(E[e^{Yx}]-1)}e^{rx}\Big)&=\sum_{n=0}^{\infty}\frac{t^{n}}{n!}D_{x}^{n}\sum_{j=0}^{\infty}\phi_{j,r}^{Y}(y)\frac{x^{j}}{j!}\label{36} \\
	&=\sum_{n=0}^{\infty}\sum_{j=0}^{\infty}\phi_{j+n,r}^{Y}(y)\frac{t^{n}}{n!}\frac{x^{j}}{j!}. \nonumber
\end{align}
Therefore, by \eqref{35} and \eqref{36}, we obtain the following theorem.
\begin{theorem}
For $n,j\ge 0$, we have
\begin{displaymath}
\phi_{j+n,r}^{Y}(y)=\sum_{l=0}^{n}\binom{n}{l}r^{n-l}\sum_{k=0}^{l}\sum_{m=0}^{j}\frac{y^{k}}{k!}\binom{j}{m}\sum_{l_{1}+\cdots+l_{k}=l}\binom{l}{l_{1},\dots,l_{k}}E\Big[S_{k}^{j-m}\prod_{i=1}^{k}Y_{i}^{l_{i}}\Big]\phi_{m,r}^{Y}(y),
\end{displaymath}
where $l_{1},l_{2},\dots l_{k}$ are positive integers.
\end{theorem}
When $Y=1$, we have
\begin{align*}
\phi_{j+n,r}(y)&=\sum_{l=0}^{n}\sum_{m=0}^{j}\binom{n}{l}\binom{j}{m}r^{n-l}\sum_{k=0}^{l}\frac{y^{k}}{k!}\sum_{l_{1}+\cdots+l_{k}=l}\binom{n}{l_{1},\dots,l_{k}}k^{j-m}\phi_{m,r}(y) \\
&=\sum_{l=0}^{n}\sum_{m=0}^{j}\sum_{k=0}^{l}\binom{n}{l}\binom{j}{m}r^{n-l}\phi_{m,r}(y)y^{k}k^{j-m}{l\brace k},
\end{align*}
where $j,n$ are nonnegative integers. \par
In particular, for $n=1$, we have
\begin{equation*}
	\phi_{j+1,r}^{Y}(y)=y\sum_{k=0}^{j}\binom{j}{k}E\big[Y^{k+1}\big]\phi_{j-k,r}^{Y}(y)+r\phi_{j,r}^{Y}(y).
\end{equation*}

\section{Conclusion}
Let $Y$ be a random variable such that the moment generating function of $Y$ exists in a neighborhood of the origin. We introduced probabilistic extensions of the $r$-Stirling numbers of the second and the $r$-Bell polynomials, namely the probabilistic $r$-Stirling numbers of the second associated with $Y$, ${n+r \brace k+r}_{r,Y}$ and the probabilistic $r$-Bell polynomials associated with $Y$, $\phi_{n,r}^{Y}(y)$. Here the latter is a natural polynomial extension of the former. Spivey found a recurrence relation for the Bell numbers $\phi_{n}$. In this paper, we generalized the recurrence relation for $\phi_{n}=\phi_{n,0}^{1}(1)$ to that for $\phi_{n,r}^{Y}(y)$ by utilizing the simple fact \eqref{15}.\par
In more detail, we obtained for ${n+r \brace k+r}_{r,Y}$ an explicit expression in terms of the $n$th moments of $S_{j}+r,\,\,(0 \le j \le k)$ in Theorem 2.1. We derived for  $\phi_{n,r}^{Y}(x)$ the generating function, a finite sum expression involving the incomplete Bell polynomials and a recurrence relation, respectively in Theorem 2.2, Theorem 2.3 and Theorem 2.4. Then the relation for the Bell numbers in \eqref{9} was generalized to first that for $\phi_{n}^{Y}=\phi_{n,0}^{Y}(1)$, second that for $\phi_{n}^{Y}(y)=\phi_{n,0}^{Y}(y)$, and finally that for $\phi_{n,r}^{Y}(y)$, respectively in Theorem 2.5, Theorem 2.6 and Theorem 2.7. \par
As one of our future projects, we would like to continue to study probabilistic extensions of many special polynomials and numbers and to find their applications to physics, science and engineering as well as to mathematics.

\end{document}